%
%

\magnification=1200

\font \fr = eufm10


\font\AAA=cmr14 at 12pt
\font\BBB=cmr14 at 8pt

\overfullrule=0in

\def\boxit#1{\hbox{\vrule
 \vtop{%
  \vbox{\hrule\kern 2pt %
     \hbox{\kern 2pt #1\kern 2pt}}%
   \kern 2pt \hrule }%
  \vrule}}

\def\half{\hbox{${1\over 2}$}}



\def\Theorem#1{\medskip\noindent {\AAA T\BBB HEOREM \rm #1.}}
\def\Prop#1{\medskip\noindent {\AAA P\BBB ROPOSITION \rm  #1.}}
\def\Cor#1{\medskip\noindent {\AAA C\BBB OROLLARY \rm #1.}}

\def\Note#1{\medskip\noindent {\AAA N\BBB OTE.\rm  #1}}

\def\Conj#1{\medskip\noindent {\AAA C\BBB ONJECTURE \rm    #1.}}

\def\Qu#1{\medskip\noindent {\AAA Q\BBB UESTION \rm    #1.}}

\def\pf{\medskip\noindent {\bf Proof.}\ }
\def\qed{\hfill  $\vrule width5pt height5pt depth0pt$}

   \def\cp{{\cal P}}

\def\cp{{\cal P}}

\def\wt{\widetilde}
\def\wh{\widehat}

\def\and{\qquad {\rm and} \qquad}

\def\bbc{{\bf C}}

\def\a{\alpha}

\def\e{\epsilon}
\def\f{\phi}
\def\g{\gamma}

\def\s{\sigma}
\def\x{\xi}
\def\z{\zeta}

\def\D{\Delta}

\def\G{\Gamma}

\def\M{\hbox{{\fr M}}}
\def\Ao{\hbox{{\fr A}}_0}
\def\f{\varphi}
\def\z{\zeta}
 \def\F{\Phi}
 \def\zb{\bar \z}

\  \vskip .2in
\centerline{\bf   RUDIN'S THEOREM AND PROJECTIVE HULLS  }
 
\vskip .2in
\centerline{\bf John Wermer}

\vskip .3in
  We denote by $\D$ the closed unit disk, by $\G$ the unit circle, and by  
 $\Ao$ the disk algebra, which consists of all functions holomorphic in
  int$(\D)$ and continuous on $\D$.  
  
  By a {\sl module over $\Ao$} we mean a vector space $\M$  of continuous complex-valued functions on $\D$ such that the constant 1 lies in $\M$, and for every $a_0\in\Ao$ and $\f\in \M$,
  one has $a_0\cdot \f\in \M$.
  
  In his book ``Real and Complex Analysis'' (1966) Walter Rudin proved Theorem 12.13 which gives, in particular:
  
  \Theorem{1} {\sl Let $\M$ be a module over $\Ao$.  Assume that the maximum principle holds for $\M$ in the following sense:
 $$
 \forall f\in\M \ {\sl and\ }\forall z_0\in {\rm int} \D,\ \ \  |f(z_0)|\leq \sup_{\G} |f|.
 \eqno{(1)}
 $$
 Then every function in $\M$ is holomorphic.}
 
 \Note{}  In 1953 Rudin had proved the corresponding statement with ``module'' replaced by ``algebra''.
 \medskip
 
 Recently, in [2], R. Harvey and H. B. Lawson introduced the notion of the {\sl projective 
 hull} of a compact set $X$ in $\bbc^n$ as follows.  Denote by $\cp_d$ the space of polynomials on $\bbc^n$ of degree $\leq d$.  Consider points $x \in \bbc^n$ such that there exists a constant $C_x$ for which 
 $$
 |P(x)|\ \leq\ (C_x)^d\, \sup_X |P|
 \eqno{(2)}$$
 for all $P\in \cp_d$ and all $d$.
 The set of such points $x$ is denoted $\wh X$ and is called the {\sl projective hull of $X$ in $\bbc^n$}.
 Note that if $C_x=1$, then (2) implies that $x$ lies in the polynomial hull of $X$.
 They make the following  
 \Conj{}\ {\sl If $\g$ is a real-analytic closed curve in $\bbc^n$, then $\wh \g - \g$ is a 1-dimensional complex analytic subvariety of $\bbc^n- \g$.}
 \Note{} The real analyticity of $\g$ is important, since the analogous conjecture for a $C^\infty$-curve fails if the curve is not pluripolar.

 In what follows we consider some questions related to Rudin's Theorem and to the projective hull of a closed curve in $\bbc^2$.
 
 \Theorem{2}  {\sl  Let $\f$ be a continuous function on $\D$
such that the restriction of $\f$ to $\G$ is real analytic.  Let
$$
\M \ =\ \{a+b\f : a,b \in \Ao\}\ \subset \ C(\D).
\eqno{(3)} 
 $$
 Assume for each $x\in {\rm int}\D$ that there exists a constant $M_x$ such that for all
 $f\in\M$:
$$
 |f(x)| \ \leq \ M_x\cdot \sup_\G |f|
 \eqno{(4)} 
 $$ 
 Then $\f$ (and hence every function in $\M$) is holomorphic in $\D$.
 }
 
 \Note{} Condition (4) states that the evaluation functional: $f\mapsto f(x)$, on $\M$, is a bounded linear functional in the sup-norm on $\G$, for all $x\in $ int$(\D)$.  It implies that two functions in $\M$ which coincide on $\G$ are equal in $\M$.  This allows us to give $\M$
 the norm: $\|f\|=\sup_\G |f|$, making $\M$ a   subspace of $C(\G)$.
 
 \medskip
 \noindent
 {\bf Proof of Theorem 2.}  Let $L$ be the functional $f\mapsto f(0)$ on $\M$.
 By (4), $\|L\|\leq M_0$, so by Hahn-Banach and F. Riesz, there exists a measure $\s$ on $\G$ such that
 $$
 f(0)\ =\ \int_\G f\,d\s\qquad\forall f\in\M.
  \eqno{(5)} 
 $$
 We do not know that $\|\s\|=1$, only that $\|\s\|<\infty$. 
 
 We take $f(\z)=\z^n$, $n>0$ and get from (5) that
  $$
  \int_\G \z^n\,d\s  \ =\ 0      \qquad n=1,2,...\ \  {\rm\ and}
  \eqno{(6)} 
 $$
  $$
  \int_\G 1\,d\s  \ =\ 1
  \eqno{(7)} 
 $$
 We do not know that $\s$ is a positive measure.
  
  We further have that
 $$
  \int_\G \z^n\,{d\theta\over 2\pi}  \ =\ 0      \qquad n=1,2,... \ \ {\rm\ and}
  \eqno{(8)} 
 $$
  $$
  \int_\G 1\,{d\theta\over 2\pi}  \ =\ 1  
   \eqno{(9)} 
 $$
 The last four equations give
  $$
  \int_\G \z^n\,   \left(d\s -     {d\theta\over 2\pi} \right)  \ =\ 0      \qquad n \ \geq \ 0.
    \eqno{(10)} 
 $$
By the F. and M. Riesz Theorem which identifies annihilating measures for the disk algebra, there exists $h\in H_0^1(\G)$ such that 
$$
 d\s -     {d\theta\over 2\pi}\ =\  h {d\theta\over 2\pi}
  \eqno{(11)} 
 $$
 Let us now fix $\f\in \M$ and apply (5) to the functions $z^n\f$, $n=0,1,2,...$ all of which belong to $\M$.  We get
 $$
  \int_\G \z^n \f \,  d\s\ =\ 0 \qquad n=1,2,...\ \ {\rm and }
 $$
  $$
  \int_\G 1\f \,  d\s\ =\ \a \ \ \ {\rm for \ some\ constant\ }\a.
 $$
 It follows that 
 $$
 \int_\G \z^n\left(\f \,d\s - \a{d\theta\over 2\pi}\right)\ =\ 0  \qquad n=0,1,2,... .
 $$
 Hence there exists $k\in H_0^1(\G)$ such that
 $$
 \f d\s -    \a {d\theta\over 2\pi}\ =\  k {d\theta\over 2\pi}
  \eqno{(12)} 
 $$
 Multiplying (11) by $\f$, we get
 $$
\f d\s -    \f {d\theta\over 2\pi}\ =\  \f h {d\theta\over 2\pi}.
  \eqno{(13)} 
 $$
 Now $h$ and $k$ are boundary functions defined on $\G$ of functions defined and holomorphic on int$\D$. (We again denote these functions on 
int$\D$ by $h$ and $k$.)  From (12) and (13) we get
$$
(\a+k){d\theta\over 2\pi}\ =\ \f(1+h){d\theta\over 2\pi}, \ \ \ {\rm as\ measures\ on\ \ } \G.
$$
 If $1+h=0$ a.e. on $\G$, then $\a+k=0$ a.e. on $\G$, so from (12) we get that
 $\f=0$ a.e. on $\G$, and so $\f$ is identically zero, which we may exclude.
 An $H^1$-function which is not identically zero is $\neq 0$ a.e. on $\G$, and so
 $1+h\neq0$ a.e. on $\G$.  It follows that
 $$
 \f\ =\ {\a+k\over 1+h} \ \ \ {\rm a.e.\ on\ \ } \G.
  \eqno{(14)} 
 $$
 Since $\f$ is real analytic on $\G$ there exists an analytic continuation $\wt \f$ of $\f$ from 
 $\G$ to an annulus: $1-\e<|\z|<1+\e$.  The function ${\a+k\over 1+h} $ is meromorphic in the strip: $1-\e<|\z|<1$ and has a non-tangential limit a.e. on $\G$.  By (14) this non-tangential limit
 $=\f(\z) = \lim_{z\to \z}\wt \f(z)$ a.e. on $\G$.  It follows that
  $$
{\a+k\over 1+h} \ =\ \wt \f\ \ \ {\rm in}\ \ \  1-\e<|\z|<1.
  \eqno{(15)} 
 $$
 Thus $  {\a+k\over 1+h} $ is analytic in that strip and is meromorphic in $\{|\z|<1\}$, 
 and hence has at most a finite number of poles in int$\D$.  We denote these poles
 by $z_1,...,z_n$.  Let 
 $$
 Q(\z)\ =\ \prod_{j=1}^n (\z-z_j).
  $$
 Then $Q\cdot \left({\a+k\over 1+h} \right)$ is analytic on int$\D$. Since ${\a+k\over 1+h} = \wt \f$
 on the strip $1-\e <|\z|<1$,   $Q\cdot \left({\a+k\over 1+h} \right)$is continuous up to $\G$, and coincides with $Q\cdot \f$ on $\G$.
 
 Thus  $Q\cdot \left({\a+k\over 1+h} \right)$ lies in $\Ao$ and hence in $\M$.
 For $\z_0\in $ int$\D$, we have, by (4), 
 $$
 \left|  Q \cdot\left({\a+k\over 1+h} \right) (\z_0) -(Q\f)(\z_0)\right|
 \ \leq\ M_{\z_0}\cdot
 \sup_\G   \left|  Q\cdot \left({\a+k\over 1+h} \right)  -(Q\f)\right|\ =\ 0.
 $$
 So $Q\left({\a+k\over 1+h} \right)$ and $(Q\f)$ agree at $\z_0$.  This holds for each
 $\z_0$.  Thus they coincide on $\D$.
 
 So $Q\f$ is analytic on int$\D$, and hence $\f$ is meromorphic on int$\D$.  But by hypothesis,
 $\f$ is continuous on $\D$. So $\f$ is analytic on $\D$ and we are done.\qed
 
 \medskip
 
 Using this notion of projective hull, we next consider a generalization of Theorem 2
 where the disk $\D$ is replaced by the {\sl punctured disk}
  $\D-\{0\}$.   Given a continuous function $\f$ defined on  $\D-\{0\}$, let $\g$ denote
 its graph over $\G$ and let $\Sigma$ denote its graph over $\D-\{0\}$.  Thus
 $$\g = \{(\z, \f(\z) : \z\in\G\}\subset \bbc^2  \qquad{\rm and}\qquad
  \Sigma = \{(\z, \f(\z) : \z\in\D-\{0\}\}\subset \bbc^2.$$
 What strengthening of hypothesis (4) will imply that $\f$ is analytic on $\D-\{0\}$ and has either a pole or a removable singularity at $\z=0$?
 
 \Prop{3}  {\sl If $\f$ is analytic on $\D-\{0\}$ and has either a pole or a removable singularity
 at $\{\z=0\}$, then $\Sigma$ is contained in the projective hull $\wh \g$ of $\g$.}

 \pf
 Fix $P\in\cp_d$ with $\sup_\g|P| \leq 1$.  Suppose $\f$  has a pole at $0$ of order $k$.
Then $P(\z, \f(\z))$ is meromorphic on$\D$ with only a pole at $\z=0$ of order $\leq dk$.

 Fix $\z_0\in $ int$\D$. The function $\z^{dk}P(\z, \f(\z))$ is holomorphic on $\D$.  Hence,
 $$
 \z_0^{dk}P(\z_0, \f(\z_0)) \ =\ {1\over2\pi i}  \int_\G{\z^{dk}P(\z, \f(\z))\over \z-\z_0} \ d\z
 $$
 Therefore 
 $$
\left| \z_0^{dk}P(\z_0, \f(\z_0))\right| \ \leq\ {1\over2\pi }  
{1 \over  1-|\z_0|}\int_\G  |P(\z, \f(\z))|\ d\theta \ \leq\   {1 \over  1-|\z_0|}.
 $$
 So
 $$
 \left|  P(\z_0, \f(\z_0))\right| \ \leq\ {1 \over  1-|\z_0|}  \left( {1\over |\z_0|^k  }   \right)^d.
 $$
 This holds for all $d$,  so $(\z_0, \f(\z_0)) \in \wh \g$.\qed
 
 \medskip
 \noindent
 {\bf Note:}\ This proposition also follows from results in [2].

 \Qu{4}  {\sl  Is the converse of Proposition 3 true? That is, given $\f$ continuous on the punctured disk and real analytic on the boundary,  with $\f$ and $\Sigma$ as before, does the hypothesis: 
 $\Sigma\subset \wh \g$
 imply that $\f$ is holomorphic on int$\D -\{0\}$ and has at $\{\z=0\}$ either a pole or a removable
 singularity?
 }
 \medskip
 
 Written explicitly, our hypothesis states:  If $\z_0 \in \D-\{0\}$, then there exits a constant
 $C_{\x_0}$ such that 
 $$
  \left|  P(\z_0, \f(\z_0))\right| \ \leq\ 
 \left(    C_{\z_0} \right)^{d}\cdot \sup_\G   \left|  P(\z, \f(\z))\right| \qquad {\rm for\ } \ P\in \cp_d
 $$

 The following theorem gives a partial answer to Question 4.  Recall that it is necessary
 that  the function $\f$ be real analytic when restricted to $\G$ (or at least that its graph 
 be a pluripolar curve in $\bbc^2$). It is natural therefore to assume $\f$ to be real
 analytic on  punctured disk. We shall  assume further that $\f$ is real analytic
 on the entire plane.
 
 \Theorem{5}  {\sl  Let $\F$ be an entire function on $\bbc^2$ written as 
 $
 \F(z,w) =\sum_{n,m=0}^\infty a_{nm} z^n w^m, 
 $
 the series converging on all of $\bbc^2$.  Let
 $$
 \f(\z) \ \equiv \ \F(\z,\zb) \ =\ \sum_{n,m=0}^\infty a_{nm} \z^n {\zb}^m,  \qquad \z\in\bbc.
 \eqno{(16)}
 $$
 Define $\g=\{(\z,\f(\z)) : \z\in \G\}$ and 
 $\Sigma=\{(\z,\f(\z)) : \z\in \D\}$.
 
 Assume that $\Sigma\subset \wh \g$.  Then $\f$ is complex analytic on $\D-\{0\}$ with a removable singularity at  $\z=0$.
 }
 
 \pf
 Suppose $\f$ is not complex analytic on  $\D-\{0\}$.  Since the series representing $\f$ converges for all $\z\in\bbc$, we have for each $R>0$ a constant $C_R$ such that
 $|a_{nm}|\leq C_R/R^{n+m}$ for all $n,m$.  Fix $d$.  We may write
 $$
 \f(\z)\ =\ \sum_{n+m\leq d} a_{nm} \z^n {\zb}^m + \sum_{n+m > d} a_{nm} \z^n {\zb}^m
 \qquad \z\in\bbc.
 $$
 We denote the second  term on the right hand side  by $\e_d(\z)$.
 
 We have for each $|\z|\leq 2$, and for each $R>2$, that
 $$\eqalign{
 |\e_d(\z)|\ &\leq \ \sum_{n+m > d}| a_{nm}| 2^{n+m}   \cr
 &\leq \ \sum_{n+m > d} {C_R\over R^{n+m}} 2^{n+m}    
\ =  \ \sum_{k > d}\sum_{n+m =k} {C_R\over R^{k}} 2^{k}   
\ =  \ \sum_{k > d}\sum_{n+m =k}  C_R(k+1)\left( {2\over R}\right)^k   \cr
&\leq\ C_R\left[ (d+1)\left( {2\over R}\right)^d + (d+2)\left( {2\over R}\right)^{d+1} +\cdots\right] \cr
 &=\ C_R\left( {2\over R}\right)^d \left[ (d+1) + (d+2)\left( {2\over R}\right) 
 + (d+3)\left( {2\over R}\right)^{2} +\cdots\right] \cr
 &\leq\ C_R\left( {2\over R}\right)^d(4+2d) \qquad{\rm for\ \ } R>4.
 }
 $$
 Hence, there exist $R_0, d_0$ with
 $$
  |\e_d(\z)|\ \leq\ C_R\left( {4\over R}\right)^d \qquad \forall \z, |\z|\leq 2\ \ {\rm and\ \ }\forall R>R_0 
{\ \rm and\ }    d>d_0.
 \eqno{(17)}
 $$
Define now for all $d$ the polynomial
$$
P_d(\z,w)\ =\ \z^dw - \sum_{n+m\leq d}  a_{nm} \z^n \z^{d-m}.
$$
 Then $P_d$ lies in $\cp_{2d}$.
 
 For $(\z,w)$ in $\g$ we have $w=\f(\z)$ and $\zb={1\over \z}$, so
 $$
 P_d(\z,w) \ =\ \z^d\left[  \sum_{n+m\leq d}a_{nm} \z^n {1\over \z^m} + \e_d(\z) \right]
 - \z^d\left[  \sum_{n+m\leq d}a_{nm} \z^n {1\over \z^m} \right] \ =\  \z^d\e_d(\z) 
 $$
 for $\z\in\G$.  Hence,  
 $$
  |P_d(\z)|\ \leq\ |\e_d(\z)|\ \leq\ C_R\left( {4\over R}\right)^d  \qquad {\rm on\ \ }\g
 \eqno{(18)}
 $$
by (17).  

Now by assumption $\f$ is not complex analytic, so for some $m>0$, $a_{nm}\neq0$.
 The function $\a\mapsto \F(\a,\bar \a)-\F(\a, {1\over\a})$ is real-analytic
 in $|\a|>0$.  Suppose $\F(\a,\bar \a)-\F(\a, {1\over\a})$ is identically zero.  Then
 $\a\mapsto \F(\a,\bar\a)$ is complex-analytic on $|\a|>0$, and so $\f$ is complex-analytic,
 contrary to assumption.
 
 Thus $\F(\a,\bar \a)-\F(\a, {1\over\a})$ is not identically zero.  We therefore can choose
 $\a_0$ in $\half <|\a| <1$ such that $\F(\a_0,\bar \a_0)-\F(\a_0, {1\over\a_0})  =  \tau\neq0$.
 We next estimate the value of $P_d$ at the point $ (\a_0,\f( \a_0))$ in $\Sigma$.
  $$\eqalign{
  P_d(\a_0,\f( \a_0)) \ &=\ \a_0^d\left[  \sum_{n+m\leq d} a_{nm} \a_0^n   \bar{\a_0}^m  
  +  \e_d(\a_0)\right]
  -\a_0^d\left[  \sum_{n+m\leq d} a_{nm} \a_0^n   \left({1\over \a_0}\right)^m  \right]  \cr
 &=\   \a_0^d\left[  \sum_{n+m\leq d} a_{nm} \a_0^n   \bar{\a_0}^m  
- \sum_{n+m\leq d} a_{nm} \a_0^n   \left({1\over \a_0}\right)^m\right]
+ \a_0^d  \e_d(\a_0)\cr
   }
  $$
 For large $d$, then, we have in view of (17) that
 $$\eqalign{
 | P_d(\a_0,\f( \a_0))|\ &\geq\ |\a_0|^d\left|  \F(\a_0, \bar \a_0)-\F(\a_0, {1\over\a_0})    \right|\cdot\half -  |\a_0|^d|\e_d(\a_0)|  \cr
 &=\ |\a_0|^d\left[  {|\tau|\over 2} - |\e_d(\a_0)|     \right] \ \geq\  |\a_0|^d{|\tau|\over 4}\cr
 }
 $$
 Thus we have 
 $$
  | P_d(\a_0,\f( \a_0))|\  \geq\  |\a_0|^d{|\tau|\over 4}  \qquad {\rm for \ large\ \ } d.
    \eqno{(19)}
  $$
  Also,  for each $R>4$, we have by (18) that $|P_d|\leq C_R({4\over R})^d$ on $\g$, so
  $$
  \left|  {  P_d \over  C_R({4\over R})^d} \right|  \ \leq  1  \ \ \ \ {\rm on\ \ } \g.
  $$  
 By our hypothesis, $(\a_0,\f( \a_0))\in \wh\g$, so by definition of $\wh\g$
 there exists a constant $C_{\a_0}$ such that
 $$
  \left|  {  P_d(\a_0,\f( \a_0)) \over  C_R({4\over R})^{d}} \right|  \ \leq  \ \left(C_{\a_0}\right)^{2d}
  \qquad {\rm for\ all\ \ }d\geq d_0.
 $$
 So
 $$
 | P_d(\a_0,\f( \a_0))|\ \leq\ C_R\left({4\over R}\cdot C_{\a_0}^2\right)^d
  \qquad {\rm for\ all\ \ }d\geq d_0.
 \eqno{(20)}
 $$
 We now choose $R$ so that ${4\over R}\cdot C_{\a_0} < |\a_0|$, and let $d\to \infty$.
 Then (20) and (19) yield a contradiction.  So $\f$ is complex-analytic on 
 $\D-\{0\}$ and hence on all of $\D$.\qed
 
 \Note{6}  Since $\f$ is real analytic on the entire $\z$-plane, it follows that $\f$ is an entire holomorphic function of $\z$.
 
 Theorem 9.2 of [2] yields that if $f$ is an  entire holomorphic function on $\bbc$, then for a closed curve $K$ on the graph $\Sigma_f$  of $f$, the projective hull $\wh K$ of $K$ 
 in $\bbc^2$ equals $K$ union the bounded components of  $\Sigma_f-K$.  Together with Theorem 5 just given we get :
 
 \Cor{7}  {\sl Let $\f$ be given by a series $\sum_{n,m=0}^\infty a_{nm} \z^n{\bar \z}^m$
 which converges  for all $\z\in\bbc$.  Let $\g=\{(\z, \f(\z)) :\z\in \G\}$. Then
 $\f$ is an entire holomorphic function of $\z$ if and only if $\{(\z, \f(\z)) :\z\in \D\} = \wh\g$.}

 \vskip .4in
 \centerline{\bf References}
 \bigskip
 
 \noindent [1]  W. Rudin, ``Real and Complex Analysis'', McGraw Hill, Inc., N.Y., 1966.
 
 \medskip
  \noindent [2] 
  F. R. Harvey and H. B. Lawson, Jr., {\sl Projective hulls and the projective Gelfand transformation},  Asian J. Math. {\bf 10}, no. 3 (2006), 279-319.


  \vskip.3in
  
  \hskip 2.5in
  Brown University
  
    \hskip 2.5in
  Providence, R.I.

  \end